\documentclass[11pt,a4paper]{article}

\usepackage{inputenc}
\usepackage{amsmath}
\usepackage{bm}
\usepackage{bbold}
\usepackage{amsthm}

\usepackage{hyperref}

\title{Efficient Parallel Algorithms\\ for Tandem Queueing System Simulation\thanks{Proc. 3rd Beijing Intern. Conf. on System Simulation and Scientific Computing, Oct. 17-19, 1995, Beijing, China: Delayed papers / Ed. by Xingren Wang et al. Chinese Assoc. for System Simulation, 1995. pp. 8-12.}}

\author{Sergei~M.~Ermakov\thanks{The research described in this
publication was made possible in part by Grant \#NVZ000 from the
International Science Foundation.} \  and Nikolai~K.~Krivulin\thanks{Faculty of Mathematics and Mechanics, St.~Petersburg State University, 28 Universitetsky Ave., St.~Petersburg, 198504, Russia, 
nkk@math.spbu.ru}
}

\date{}

\newcounter{algorithm}
\newenvironment{algorithm}{\bf\refstepcounter{algorithm}\par
\vspace{1.5ex}%
{\bf Algorithm~\thealgorithm.}\nopagebreak[4]
\vspace{-1.0ex}
\begin{tabbing}}{\end{tabbing}}

\begin{document}

\maketitle

\begin{abstract}
Parallel algorithms designed for simulation and performance evaluation of
single-server tandem queueing systems with both infinite and finite buffers
are presented. The algorithms exploit a simple computational procedure based
on recursive equations as a representation of system dynamics. A brief
analysis of the performance of the algorithms are given to show that they
involve low time and memory requirements.
\end{abstract}

\section{Introduction}

The simulation of a queueing system is normally an iterative process which
involves generation of random variables associated with current events in the
system, and evaluation of the system state variables when new events occur
\cite{Erma75,Chen90,Gree91}. In a system being simulated the random variables
may represent the interarrival and service time of customers, whereas, as
state variables, the arrival and departure time of customers, and the service
initiation and completion time can be considered.

The usual way to represent dynamics of queueing systems as well as their
performance criteria relies on recursive equations describing evolution of
system state variables \cite{Chen90,Gree91,Kriv93,Chen94,Kriv94}. Since the
recursive equations actually determine a global structure of changes in the
system state variables consecutively, they can serve as a basis for the
development of efficient simulation algorithms
\cite{Chen90,Gree91,Erma94,Kriv94}.

In this paper, we assume as in \cite{Gree91,Erma94,Kriv94} that appropriate
realizations of the random variables involved in simulation are available when
required, and we therefore concentrate only on deterministic parallel
algorithms of evaluating the system state variables from these realizations.
Methods and algorithms of generating random variables and their analysis
can be found in \cite{Erma75}. A thorough investigation of parallel simulation
from the viewpoint of statistics is given in \cite{Heid88}.

We present parallel algorithms designed for simulation and performance
evaluation of open single-server tandem queueing systems with both infinite
and finite buffers. The algorithms are based on a simple computational
procedure which exploits a particular order of evaluating the system state
variables from the related recursive equations, and they are intended for
implementation on either a vector processor or single instruction, multiple
data (SIMD) parallel processors \cite{Orte88}. The analysis of their
performance shows that the algorithms involve low time and memory
requirements.

In Section~\ref{s-MTQ}, we give recursive equations which describe the
dynamics of tandem systems with both infinite and finite buffers. Furthermore,
tandem system performance criteria are represented in terms of state variables
involved in the recursive equations. In Section~\ref{s-TQSA}, parallel
simulation algorithms are presented and their performance is discussed. A
brief conclusion is given in Section~\ref{s-C}.

\section{Models of Tandem Queues} \label{s-MTQ}

In this section we consider recursive equation based models of tandem queues,
and give related representation of system performance measures. We start with
a simple model of a single-server tandem queueing system with infinite
buffers, and then extend it to more complicated models of systems with finite
buffers, in which servers may be blocked according to some blocking rule.

\subsection{Tandem Queues with Infinite Buffers}

\begin{figure}[hhh]
\begin{center}
\begin{picture}(55,15)

\put(6,7){$0$}
\put(18,7){$1$}
\put(41,7){$N$}

\put(1,5){\line(1,0){4}}
\put(13,5){\line(1,0){4}}
\put(37,5){\line(1,0){4}}

\put(1,1){\line(1,0){4}}
\put(13,1){\line(1,0){4}}
\put(37,1){\line(1,0){4}}

\put(5,5){\line(0,-1){4}}
\put(17,5){\line(0,-1){4}}
\put(41,5){\line(0,-1){4}}

\put(9,3){\vector(1,0){4}}
\put(21,3){\vector(1,0){4}}
\put(33,3){\vector(1,0){4}}
\put(45,3){\vector(1,0){6}}

\put(7,3){\circle{3}}
\put(19,3){\circle{3}}
\put(43,3){\circle{3}}

\multiput(27,3)(2,0){3}{\circle*{1}}

\end{picture}

\end{center}
\caption{Tandem queues with infinite buffers.} \label{f-TQIB}
\end{figure}
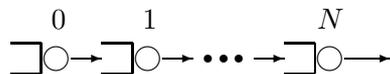

Consider a series of $\, N \,$ single-server queues with infinite buffers,
depicted in Fig.~\ref{f-TQIB}. An additional queue labelled with $\, 0 \,$ is
included in the model to represent the external arrival stream of customers.

Each customer that arrives into the system is initially placed in the buffer
at the $1$st server and then has to pass through all the queues one after the
other. Upon the completion of his service at server $\, n $, the customer is
instantaneously transferred to queue $\, n+1 $, $\, n=1,\dots,N-1 $, and
occupies the $(n+1)$st server provided that it is free. If the customer finds
this server busy, he is placed in its buffer and has to wait until the service
of all his predecessors is completed.

For each queue $\, n \,$ in the system, $\, n=0,1,\ldots, N $, we introduce
the following notations:
{\parskip=1.5ex

$ A_{n}^{k} $, \quad the $k$th arrival epoch to the queue;

$ B_{n}^{k} $, \quad the $k$th service initiation time at the queue;

$ C_{n}^{k} $, \quad the $k$th service completion time at the queue;

$ D_{n}^{k} $, \quad the $k$th departure epoch from the queue.
\par}
\vspace{1ex}

Furthermore, let us denote the time between the arrivals of $k$th customer and
his predecessor to the system by $\, \tau_{0}^{k} $, and the service time of
the $k$th customer at server $\, n \,$ by $\, \tau_{n}^{k} $,
$\, n=1,\dots,N $, $\, k=1,2,\dots $. We assume that $\, \tau_{n}^{k} \geq 0 $
are given parameters, whereas $\, A_{n}^{k}, B_{n}^{k}, C_{n}^{k} $, and
$\, D_{n}^{k} \,$ present unknown state variables. Finally,
for each $\, n=0,\dots,N $, we define $\, D_{n}^{k} \equiv 0 \,$ for all
$\, k \leq 0 $, and $\, D_{-1}^{k} \equiv 0 \,$ for all $\, k=1,2,\ldots $.

With the condition that the system starts operating at time zero, and its
servers are free of customers at the initial time, the state variables
in the model can be related by the equations
\cite{Chen90,Gree91,Chen94,Kriv94}
\begin{eqnarray*}
A_{n}^{k} & = & D_{n-1}^{k}, \\
B_{n}^{k} & = & A_{n}^{k} \vee D_{n}^{k-1}, \\
C_{n}^{k} & = & B_{n}^{k} + \tau_{n}^{k}, \\
D_{n}^{k} & = & C_{n}^{k},
\end{eqnarray*}
where the symbol $\, \vee \,$ stands for the maximum operator,
$\, n=0,1,\dots,N $, $\, k=1,2,\dots $. Clearly, the above set of recursive
equations may be reduced to two equations
\begin{eqnarray}
B_{n}^{k} & = & D_{n-1}^{k} \vee D_{n}^{k-1}, \label{e1-IB} \\
D_{n}^{k} & = & B_{n}^{k} + \tau_{n}^{k},     \label{e2-IB}
\end{eqnarray}
and even to the equation
\begin{equation} \label{e-IB1}
D_{n}^{k} = (D_{n-1}^{k} \vee D_{n}^{k-1}) + \tau_{n}^{k},
\end{equation}
which will provide the basic representations for simulation algorithms in the
next sections.

\subsection{Tandem Queues with Finite Buffers}

Suppose now that the buffers of servers in the open tandem system have finite
capacity. Furthermore, we assume that the servers may be blocked according to
some blocking rule. In this paper, we restrict our consideration to
{\it manufacturing} blocking and {\em communication} blocking which are most
commonly encountered in practice \cite{Chen90,Gree91,Chen94}.

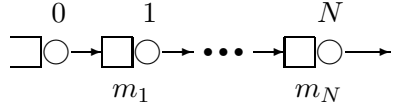
\begin{figure}[hhh]
\begin{center}
\begin{picture}(55,15)

\put(6,10){$0$}
\put(18,10){$1$}
\put(41,10){$N$}

\put(1,8){\line(1,0){4}}
\put(13,8){\line(1,0){4}}
\put(37,8){\line(1,0){4}}

\put(1,4){\line(1,0){4}}
\put(13,4){\line(1,0){4}}
\put(37,4){\line(1,0){4}}

\put(5,8){\line(0,-1){4}}
\put(13,8){\line(0,-1){4}}
\put(17,8){\line(0,-1){4}}
\put(37,8){\line(0,-1){4}}
\put(41,8){\line(0,-1){4}}

\put(9,6){\vector(1,0){4}}
\put(21,6){\vector(1,0){4}}
\put(33,6){\vector(1,0){4}}
\put(45,6){\vector(1,0){6}}

\put(7,6){\circle{3}}
\put(19,6){\circle{3}}
\put(43,6){\circle{3}}

\multiput(27,6)(2,0){3}{\circle*{1}}

\put(14,0){$m_{1}$}
\put(38,0){$m_{N}$}

\end{picture}
\end{center}
\caption{Tandem queues with finite buffers.}\label{f-TQFB}
\end{figure}

Let us consider an open tandem system of $\, N \,$ queues (Fig.~\ref{f-TQFB}),
and assume the buffer at the $n$th server, $\, n=1,\ldots,N $, to be of the
capacity $\, m_{n} $, $\, 0 < m_{n} < \infty $.

{\it Manufacturing Blocking.} First we suppose that the dynamics of the system
follows the manufacturing blocking rule. Under this type of blocking, if upon
completion of a service, the $n$th server sees the buffer of the $(n+1)$st
server full, it cannot be unoccupied and has to be busy until the $(n+1)$st
server completes its current service to provide a free space in its buffer.
Clearly, since the customers leave the system upon their service completion at
the $N$th server, this server cannot be blocked.

With the additional condition that $\, D_{n}^{k} \equiv 0 \,$ if $\, n > N $,
one can describe the dynamics of the system by the equations
\cite{Chen90,Gree91,Chen94,Kriv94}
\begin{eqnarray}
B_{n}^{k} & = & D_{n-1}^{k} \vee D_{n}^{k-1},          \label{e1-FBM} \\
C_{n}^{k} & = & B_{n}^{k} + \tau_{n}^{k},              \label{e2-FBM} \\
D_{n}^{k} & = & C_{n}^{k} \vee D_{n+1}^{k-m_{n+1}-1}.  \label{e3-FBM}
\end{eqnarray}

{\it Communication Blocking.} This rule does not permit a server to initiate
service of a customer if the buffer of the next server is full. In that case,
the server remains unavailable until the current service at the next server is
completed.

Let us assume that the system depicted in Fig.~\ref{f-TQFB} follows
communication blocking, and introduce the notation $\, H_{n}^{k} \,$ to denote
the time instant at which the $n$th server becomes ready to check whether
there is empty space at the buffer of the next server, and to initiate service
of customer $\, k \,$ if it is possible. Now the system dynamics may be
represented by the equations \cite{Chen90,Chen94,Kriv94}
\begin{eqnarray}
H_{n}^{k} & = & D_{n-1}^{k} \vee D_{n}^{k-1},          \label{e1-FBC} \\
B_{n}^{k} & = & H_{n}^{k} \vee D_{n+1}^{k-m_{n+1}-1},  \label{e2-FBC} \\
D_{n}^{k} & = & B_{n}^{k} + \tau_{n}^{k}.              \label{e3-FBC}
\end{eqnarray}

\subsection{Representation of System Performance}

Suppose that we observe the system until the $K$th service completion at
server $\, n $, $\, 1 \leq n \leq N $. As is customary in queueing system
simulation, we assume that $\, K > N $. The following average quantities are
normally considered as performance criteria for server $\, n \,$ in the
observation period \cite{Chen90,Kriv90,Kriv93,Kriv94}:
\begin{tabbing}
\qquad \= system time \\
       \> of one customer: \quad
\= $ S_{n} = \sum_{k=1}^{K} (D_{n}^{k}-A_{n}^{k})/K $, \\ \\
       \> waiting time \\
       \> of one customer:
\> $ W_{n} = \sum_{k=1}^{K} (B_{n}^{k}-A_{n}^{k})/K $, \\ \\
       \> throughput rate \\
       \> of the server:
\> $ T_{n} = K/D_{n}^{K} $, \\ \\
       \> utilization \\
       \> of the server:
\> $ U_{n} = \sum_{k=1}^{K} \tau_{n}^{k}/D_{n}^{K} $, \\ \\
       \> number of \\
       \> customers:
\> $ J_{n} = \sum_{k=1}^{K} (D_{n}^{k}-A_{n}^{k})/D_{n}^{K} $, \\ \\
       \> queue length \\
       \> at the server:
\> $ Q_{n} = \sum_{k=1}^{K} (B_{n}^{k}-A_{n}^{k})/D_{n}^{K} $.
\end{tabbing}

Clearly, the above criteria are suited to the systems with both infinite and
finite buffers. Furthermore, one can consider the average idle time of server
$\, n $, which presents a criterion inherent only in the systems with finite
buffers. It is defined for the manufacturing and communication blocking rules
respectively as \cite{Chen90,Kriv94}
\begin{eqnarray*}
IM_{n} & = & \sum_{k=1}^{K} (D_{n}^{k}-C_{n}^{k})/K, \\
IC_{n} & = & \sum_{k=1}^{K} (B_{n}^{k}-H_{n}^{k})/K.
\end{eqnarray*}
Note finally that these expressions may be also written in terms of departure
epochs and service times in the same form as
$$
I_{n} = \sum_{k=1}^{K} \left(D_{n}^{k}-(D_{n-1}^{k} \vee D_{n}^{k-1})
      - \tau_{n}^{k}\right)/K.
$$

\section{Tandem Queues Simulation Algorithms} \label{s-TQSA}

We start with the description of a simple simulation procedure designed for
the tandem system with infinite buffers, and then extend the procedure to
algorithms for systems with finite buffers and blocking. It is shown how the
algorithms can be refined so as to evaluate system performance. In addition,
time and memory requirements associated with the algorithms are briefly
discussed.

\subsection{The Basic Simulation Procedure}

We use the procedure proposed in \cite{Erma94}, which was designed for the 
simulation of the tandem queueing system described by equations 
(\ref{e1-IB}-\ref{e2-IB}). It actually performs computations of successive 
state variables $\, B_{n}^{k} \,$ and $\, D_{n}^{k} \,$ with indices being 
varied in a particular order. According to this order, at each iteration 
$\, i $, the variables with $\, n+k=i $, $ i=1,2,\ldots $, have to be 
evaluated. The next algorithm shows how to implement this procedure to the 
simulation of the first $\, K \,$ customers in a tandem queueing system with
infinite buffers and $\, N \,$ servers, $\, K > N $.
\begin{algorithm} \label{a-IBS}%
Set \= $\, d_{i} = 0, \;\; i=-1,0,\ldots,N; $ \\
for \> $\, i=1,\ldots,K+N, $ do \\
    \> $\, j_{0} \longleftarrow \max(1,i-N); $ \\
    \> $\, J \longleftarrow \min(i,K); $ \\
    \> for \= $\, j=j_{0},j_{0}+1,\ldots,J, $ do \\
    \>     \>  $\, b_{i-j} \longleftarrow d_{i-j-1} \vee d_{i-j}; $ \\
    \>     \>  $\, d_{i-j} \longleftarrow b_{i-j} + \tau_{i-j}^{j}. $
\end{algorithm}
In Algorithm~\ref{a-IBS}, the variables $\, b_{n} \,$ and $\, d_{n} \,$ serve
all the iterations to store current values of $\, B_{n}^{k} \,$ and
$\, D_{n}^{k} \,$ respectively, for $\, k=1,\ldots,K $. Upon the completion of
the algorithm, we have for server $\, n \,$ the $K$th departure time saved in
$\, d_{n} $, $\, n=0,1,\ldots,N $.

Since one maximization and one addition have to be performed so as to get new
variables $\, B_{n}^{k} \,$ and $\, D_{n}^{k} $, one can conclude that the
entire algorithm requires $\, O(2(N+1)K) \,$ arithmetic operations without
considering index manipulations. Moreover, the order in which the variables
are evaluated within each iteration is essential for reducing memory used for
computations. It is easy to see that only $\, O(N+1) \,$ memory locations are
actually required with this order, provided only the departure epochs
$\, D_{n}^{k} \,$ are to be calculated. To illuminate the memory requirements,
let us represent Algorithm~\ref{a-IBS} in another form as
\begin{algorithm} \label{a-IBS1}%
Set \= $\, d_{i} = 0, \;\; i=-1,0,\ldots,N; $ \\
for \> $\, i=1,\ldots,K+N $, do \\
    \> $\, j_{0} \longleftarrow \max(1,i-N); $ \\
    \> $\, J \longleftarrow \min(i,K); $ \\
    \> for \= $\, j=j_{0},j_{0}+1,\ldots,J, $ do \\
    \>     \> $\, d_{i-j} \longleftarrow d_{i-j-1} \vee d_{i-j}
                                                           + \tau_{i-j}^{j}. $
\end{algorithm}

Finally, we suppose that there is a computer system with either a vector
processor or SIMD parallel processors available for tandem queueing system
simulation. In that case, we can use the following algorithm, which is
actually a simple modification of Algorithm~\ref{a-IBS}.
\begin{algorithm} \label{a-IBP}%
Set \= $\, d_{i} = 0, \;\; i=-1,0,\ldots,N; $ \\
for \> $\, i=1,\ldots,K+N, $ do \\
    \> $\, j_{0} \longleftarrow \max(1,i-N); $ \\
    \> $\, J \longleftarrow \min(i,K); $ \\
    \> in parallel, for \= $\, j=j_{0},j_{0}+1,\ldots,J, $ do \\
    \>     \>  $\, b_{i-j} \longleftarrow d_{i-j-1} \vee d_{i-j}; $ \\
    \> in parallel, for \> $\, j=j_{0},j_{0}+1,\ldots,J, $ do \\
    \>     \>  $\, d_{i-j} \longleftarrow b_{i-j} + \tau_{i-j}^{j}. $
\end{algorithm}

Let $\, P \,$ denote the length of vector registers of the vector processor or
the number of parallel processors, depending on whether a vector or parallel
computer system appears to be available. It is not difficult to see that
Algorithm~\ref{a-IBP} requires the condition $\, P \geq N+1 \,$ to be
satisfied. Otherwise, if $\, P < N+1 $, one simply has to rearrange
computations so as to execute each iteration in several parallel steps. In
other words, all operations within an iteration should be sequentially
separated into groups of $\, P \,$ operations, assigned to the sequential
steps.

It has been shown in \cite{Erma94} that for any integer $\, P > 0 $,
Algorithm~\ref{a-IBP} requires
$\, O\left(2N+2K+2\lfloor N/P\rfloor(K-P)\right) \,$ parallel (vector)
operation, where $\, \lfloor x\rfloor \,$ denotes the greatest integer less
than or equal to $\, x $. Moreover, provided that $\, P = N+1 $, the algorithm
achieves linear speedup in relation to Algorithm~\ref{a-IBS} as the number of
customers $\, K \to \infty $. Finally, it is easy to understand that
Algorithm~\ref{a-IBP} entails $\, O(2(N+1)) \,$ memory locations.

\subsection{Simulation of Queues with Finite Buffers}

Taking equations (\ref{e1-FBM}-\ref{e3-FBM}) as the starting point,
we can readily rewrite Algorithm~\ref{a-IBP} so as to make it possible to
simulate tandem queueing systems with manufacturing blocking. Let us first
introduce the variables $\, b_{n} \,$ and $\, c_{n} \,$ to represent current
values of $\, B_{n}^{k} \,$ and $\, C_{n}^{k} $. Since calculation of
$\, D_{n}^{k} \,$ involves taking account of the value of
$\, D_{n+1}^{k-m_{n+1}-1} $, one has to keep in memory all values
$\, D_{n+1}^{j} \,$ with $\, j=k-m_{n+1}-1,k-m_{n+1},\ldots,k $. Therefore, we
further introduce the variables $\, d_{n}^{j} \,$ as memory locations of these
values, $\, n=1,\ldots,N $, $\, j=0,1,\ldots,m_{n} $. The locations
$\, d_{n}^{0},d_{n}^{1},\ldots,d_{n}^{m_{n}} \,$ are intended to be occupied
using cyclic overwriting so that the value $\, D_{n}^{k} \,$ is put into the
location $\, d_{n}^{j} \,$ with $\, j = k \bmod (m_{n}+1) $, where
$\, \bmod \,$ indicates the modulo operation.

In order to simplify further formulas, we define the index function
$$
\rho(k,n) = \left\{
	     \begin{array}{ll}
	       k \bmod (m_{n}+1), & \mbox{if $\, 1 \leq n \leq N $} \\
	       0,                 & \mbox{otherwise}
	     \end{array}
	    \right.
$$
for all $\, k=1,2,\ldots $. Finally, with the variables $\, d_{-1}^{0} $,
$\, d_{0}^{0} $, and $\, d_{N+1}^{0} \,$ reserved respectively for
$\, D_{-1}^{0} $, $\, D_{0}^{k} $, and $\, D_{N+1}^{k} $, we have the next
parallel algorithm.
\begin{algorithm} \label{a-FBMP}%
Set \= $\, d_{-1}^{0},d_{0}^{0},d_{N+1}^{0} = 0; $ \\
set \> $\, d_{i}^{j} = 0, \;\; i=1,\ldots,N, \;\; j=0,1,\ldots,m_{i}; $ \\
for \> $\, i=1,\ldots,K+N, $ do \\
    \> $\, j_{0} \longleftarrow \max(1,i-N); $ \\
    \> $\, J \longleftarrow \min(i,K); $ \\
    \> in parallel, for \= $\, j=j_{0},j_{0}+1,\ldots,J, $ do \\
    \>     \>  $\, b_{i-j} \longleftarrow d_{i-j-1}^{\rho(j,i-j-1)}
                               \vee d_{i-j}^{\rho(j-1,i-j)}; $ \\
    \> in parallel, for \> $\, j=j_{0},j_{0}+1,\ldots,J, $ do \\
    \>     \>  $\, c_{i-j} \longleftarrow b_{i-j} + \tau_{i-j}^{j}; $ \\
    \> in parallel, for \> $\, j=j_{0},j_{0}+1,\ldots,J, $ do \\
    \>     \>  $\, d_{i-j}^{\rho(j,i-j)} \longleftarrow c_{i-j}
                            \vee d_{i-j+1}^{\rho(j,i-j+1)}. $
\end{algorithm}

Consider now equations (\ref{e1-FBC}-\ref{e3-FBC}) which describe the
dynamics of the tandem system operating under the communication blocking rule.
With the variables $\, h_{n} $, $\, n=0,1,\ldots,N $, used as storage for the
values of $\, H_{n}^{k} $, it is easy to arrive at
\begin{algorithm} \label{a-FBCP}%
Set \= $\, d_{-1}^{0},d_{0}^{0},d_{N+1}^{0} = 0; $ \\
set \> $\, d_{i}^{j} = 0, \;\; i=1,\ldots,N, \;\; j=0,1,\ldots,m_{i}; $ \\
for \> $\, i=1,\ldots,K+N, $ do \\
    \> $\, j_{0} \longleftarrow \max(1,i-N); $ \\
    \> $\, J \longleftarrow \min(i,K); $ \\
    \> in parallel, for \= $\, j=j_{0},j_{0}+1,\ldots,J, $ do \\
    \>     \>  $\, h_{i-j} \longleftarrow d_{i-j-1}^{\rho(j,i-j-1)}
                               \vee d_{i-j}^{\rho(j-1,i-j)}; $ \\
    \> in parallel, for \> $\, j=j_{0},j_{0}+1,\ldots,J, $ do \\
    \>     \>  $\, b_{i-j} \longleftarrow h_{i-j}
                               \vee d_{i-j+1}^{\rho(j,i-j+1)}; $ \\
    \> in parallel, for \> $\, j=j_{0},j_{0}+1,\ldots,J, $ do \\
    \>     \>  $\, d_{i-j}^{\rho(j,i-j)} \longleftarrow b_{i-j}
                                         + \tau_{i-j}^{j}. $
\end{algorithm}

In fact, both algorithms differ from Algorithm~\ref{a-IBP} in that at every
iteration, they involve three parallel operations each, whereas the latter
does two operations. Therefore, we may extend the above estimate of time
requirements for Algorithm~\ref{a-IBP} to Algorithm~\ref{a-FBMP} and
Algorithm~\ref{a-FBCP}, which then becomes
$\, O\left(3N+3K+3\lfloor N/P\rfloor(K-P)\right) $. The number of memory
locations now involved in computations can be evaluated as
$\, O(3(N+1)+M+1) $, where $\, M = \sum_{i=1}^{N} m_{i} $.

\subsection{Evaluation of Performance Criteria}

In order to present a modification of Algorithm~\ref{a-IBP} suitable for the
evaluation of the tandem system performance criteria introduced in the
previous section, first define the additional variables $\, x_{n} $,
$\, y_{n} $, and $\, z_{n} $, $\, n=0,1,\ldots,N $, to represent the memory
locations which are to store current values of the sums
$$
\sum_{i=1}^{k} (A_{n}^{i}-D_{n}^{i}), \qquad
\sum_{i=1}^{k} (B_{n}^{i}-A_{n}^{i}), \qquad
\sum_{i=1}^{k} \tau_{n}^{i},
$$
respectively. Taking into account that in the tandem systems with both
infinite and finite buffers, we have $\, A_{n}^{k} = D_{n-1}^{k} \,$ for all
$\, n=0,1,\ldots,N $, and $\, k=1,2,\ldots $, we may write the following
parallel algorithm.
\begin{algorithm} \label{a-IBPP}%
Set \= $\, d_{i} = 0, \;\; i=-1,0,\ldots,N; $ \\
set \> $\, x_{i}, y_{i}, z_{i} = 0, \;\; i=0,1,\ldots,N; $ \\
for \> $\, i=1,\ldots,K+N, $ do \\
    \> $\, j_{0} \longleftarrow \max(1,i-N); $ \\
    \> $\, J \longleftarrow \min(i,K); $ \\
    \> in parallel, for \= $\, j=j_{0},j_{0}+1,\ldots,J, $ do \\
    \>     \>  $\, x_{i-j} \longleftarrow x_{i-j} - d_{i-j-1}; $ \\
    \> in parallel, for \> $\, j=j_{0},j_{0}+1,\ldots,J, $ do \\
    \>     \>  $\, y_{i-j} \longleftarrow y_{i-j} - d_{i-j-1}; $ \\
    \> in parallel, for \> $\, j=j_{0},j_{0}+1,\ldots,J, $ do \\
    \>     \>  $\, b_{i-j} \longleftarrow d_{i-j-1} \vee d_{i-j}; $ \\
    \> in parallel, for \> $\, j=j_{0},j_{0}+1,\ldots,J, $ do \\
    \>     \>  $\, d_{i-j} \longleftarrow b_{i-j} + \tau_{i-j}^{j}; $ \\
    \> in parallel, for \> $\, j=j_{0},j_{0}+1,\ldots,J, $ do \\
    \>     \>  $\, x_{i-j} \longleftarrow x_{i-j} + d_{i-j}; $ \\
    \> in parallel, for \> $\, j=j_{0},j_{0}+1,\ldots,J, $ do \\
    \>     \>  $\, y_{i-j} \longleftarrow y_{i-j} + b_{i-j}; $ \\
    \> in parallel, for \> $\, j=j_{0},j_{0}+1,\ldots,J, $ do \\
    \>     \>  $\, z_{i-j} \longleftarrow z_{i-j} + \tau_{i-j}^{j}. $
\end{algorithm}

As it is easy to understand, Algorithm~\ref{a-IBPP} requires
$\, O\left(7N+7K+7\lfloor N/P\rfloor(K-P)\right) \,$ parallel operations, and
involves $\, O(5(N+1)) \,$ memory locations. Upon the completion of the
algorithm, the performance criteria associated with each queue $\, n $,
$\, n=1,\ldots,N $, can be calculated as
\begin{eqnarray*}
& S_{n} = x_{n}/K, \qquad & W_{n} = y_{n}/K, \\
& T_{n} = K/d_{n}, \qquad & U_{n} = x_{n}/d_{n}, \\
& J_{n} = x_{n}/d_{n}, \qquad & Q_{n} = y_{n}/d_{n}.
\end{eqnarray*}

One can modify both Algorithm~\ref{a-FBMP} and Algorithm~\ref{a-FBCP} to
provide performance evaluation in tandem queueing systems with finite buffers
in an analogous way. Specifically, the next two algorithms intended to compute
the average idle time of each server in the system. The first one based on
Algorithm~\ref{a-FBMP} is designed for the system operating under the
manufacturing blocking rule.
\begin{algorithm} \label{a-FBMPP}%
Set \= $\, d_{-1}^{0},d_{0}^{0},d_{N+1}^{0} = 0; $ \\
set \> $\, d_{i}^{j} = 0, \;\; i=1,\ldots,N, \;\; j=0,1,\ldots,m_{i}; $ \\
set \> $\, x_{i} = 0, \;\; i=0,1,\ldots,N; $ \\
for \> $\, i=1,\ldots,K+N, $ do \\
    \> $\, j_{0} \longleftarrow \max(1,i-N); $ \\
    \> $\, J \longleftarrow \min(i,K); $ \\
    \> in parallel, for \= $\, j=j_{0},j_{0}+1,\ldots,J, $ do \\
    \>     \>  $\, b_{i-j} \longleftarrow d_{i-j-1}^{\rho(j,i-j-1)}
                               \vee d_{i-j}^{\rho(j-1,i-j)}; $ \\
    \> in parallel, for \> $\, j=j_{0},j_{0}+1,\ldots,J, $ do \\
    \>     \>  $\, c_{i-j} \longleftarrow b_{i-j} + \tau_{i-j}^{j}; $ \\
    \> in parallel, for \> $\, j=j_{0},j_{0}+1,\ldots,J, $ do \\
    \>     \>  $\, x_{i-j} \longleftarrow x_{i-j} - c_{i-j}; $ \\
    \> in parallel, for \> $\, j=j_{0},j_{0}+1,\ldots,J, $ do \\
    \>     \>  $\, d_{i-j}^{\rho(j,i-j)} \longleftarrow c_{i-j}
                            \vee d_{i-j+1}^{\rho(j,i-j+1)}; $ \\
    \> in parallel, for \> $\, j=j_{0},j_{0}+1,\ldots,J, $ do \\
    \>     \>  $\, x_{i-j} \longleftarrow x_{i-j} + d_{i-j}^{\rho(j,i-j)}. $
\end{algorithm}
The variable $\, x_{n} \,$ inserted in Algorithm~\ref{a-FBMPP} serves for
each $\, n $, $\, n=0,1,\ldots,N \,$ to represent current values of the sums
$\, \sum_{i=1}^{k} (D_{n}^{i}-C_{n}^{i}) $. Upon the completion of the
algorithm, one can calculate $\, x_{n}/K \,$ which gives the value of
$\, IM_{n} $. The time and memory costs can be estimated respectively as
$\, O\left(5N+5K+5\lfloor N/P\rfloor(K-P)\right) \,$ and
$\, O(4(N+1)+M+1) $.

\begin{algorithm} \label{a-FBCPP}%
Set \= $\, d_{-1}^{0},d_{0}^{0},d_{N+1}^{0} = 0; $ \\
set \> $\, d_{i}^{j} = 0, \;\; i=1,\ldots,N, \;\; j=0,1,\ldots,m_{i}; $ \\
set \> $\, x_{i} = 0, \;\; i=0,1,\ldots,N; $ \\
for \> $\, i=1,\ldots,K+N, $ do \\
    \> $\, j_{0} \longleftarrow \max(1,i-N); $ \\
    \> $\, J \longleftarrow \min(i,K); $ \\
    \> in parallel, for \= $\, j=j_{0},j_{0}+1,\ldots,J, $ do \\
    \>     \>  $\, h_{i-j} \longleftarrow d_{i-j-1}^{\rho(j,i-j-1)}
                               \vee d_{i-j}^{\rho(j-1,i-j)}; $ \\
    \> in parallel, for \> $\, j=j_{0},j_{0}+1,\ldots,J, $ do \\
    \>     \>  $\, x_{i-j} \longleftarrow x_{i-j} - h_{i-j}; $ \\
    \> in parallel, for \> $\, j=j_{0},j_{0}+1,\ldots,J, $ do \\
    \>     \>  $\, b_{i-j} \longleftarrow h_{i-j}
                               \vee d_{i-j+1}^{\rho(j,i-j+1)}; $ \\
    \> in parallel, for \> $\, j=j_{0},j_{0}+1,\ldots,J, $ do \\
    \>     \>  $\, x_{i-j} \longleftarrow x_{i-j} + b_{i-j}; $ \\
    \> in parallel, for \> $\, j=j_{0},j_{0}+1,\ldots,J, $ do \\
    \>     \>  $\, d_{i-j}^{\rho(j,i-j)} \longleftarrow b_{i-j}
                                         + \tau_{i-j}^{j}. $
\end{algorithm}
With the same time and memory requirements as for the previous algorithm,
Algorithm~\ref{a-FBCPP} allows one to evaluate the average idle time of each
servers in tandem queues with communication blocking. It produces the sums
$\, \sum_{i=1}^{K} (B_{n}^{i}-H_{n}^{i}) \,$ stored in $\, x_{n} $,
$\, n=1,\ldots,N $, which can be used in calculation of the criteria
$\, IC_{n} \,$ with the expression $\, x_{n}/K $.

\section{Conclusions} \label{s-C}

Parallel algorithms which offer a quite simple and efficient way of simulating
tandem queueing system have been proposed. It has been shown that the
algorithms involve low time and memory requirements. Specifically, one can
conclude that the parallel simulation of the first $\, K \,$ customers in a
system with $\, N \,$ queues requires the time of order
$\, O\left(L(N+K+\lfloor N/P\rfloor(K-P))\right) $, where $\, P \,$ is the
number of processors, $\, L \,$ is a small constant comparable with the number
of the performance criteria being evaluated. Note, however, that this estimate
ignores the time required for computing indices, and allocating and moving
data, which can have an appreciable effect on the performance of parallel
algorithms in practice.

\bibliographystyle{utphys}

\bibliography{Efficient_parallel_algorithms_for_tandem_queueing_system_simulation}

\end{document}